*Посвящается академику А.А. Петрову (1934 – 2011)*

# Заметка об эффективной вычислимости конкурентных равновесий в транспортно-экономических моделях

*Гасников А.В.*

gasnikov@yandex.ru

Кафедра математических основ управления,

Лаборатория структурных методов анализа данных в предсказательном моделировании,

Факультет управления и прикладной математики МФТИ

**Аннотация**

В данной заметке мы рассмотрим определенный класс транспортно-экономических задач, в которых поиск конкурентного равновесия может быть эффективно реализован, поскольку сводится к поиску седловой точки в выпукло-вогнутой игре. Важно отметить наличие эволюционной интерпретации возникающего конкурентного равновесия.

**Ключевые слова:** конкурентное равновесие, эволюционная игра, седловая точка, многоуровневая оптимизация, равновесие макросистемы.

## 1. Введение

В данной статье мы сосредоточим внимание на транспортно-экономических моделях, объединяющих в себе, в частности, модели из недавних работ [1], [2]. Работа мотивирована обоснованием существующих и созданием новых моделей транспортного планирования, включающих модели роста транспортной инфраструктуры городов, формирования матрицы корреспонденций и равновесного распределения потоков.

Имеется ориентированный транспортный граф, каждое ребро которого характеризуется неубывающей функцией затрат $\tau_e(f_e)$ на прохождение этого ребра, в зависимости от потока по этому ребру. Можно еще ввести затраты на прохождения вершин графа $E$, но это ничего не добавляет с точки зрения последующих математических выкладок [1]. Часть вершин графа является источниками, часть стоками (эти множества вершин могут пересекаться). В источниках $O$ и стоках $D$ имеются (соответственно) пункты производства и пункты потребления. Для большей наглядности в первой половине статьи мы будем считать, что производится и потребляется лишь один продукт. Несложно все, что далее будет написано, обобщить на многопродуктовый рынок.

Задача разбивается на две подзадачи разного уровня [3]. На нижнем уровне, соответствующем быстрому времени, при заданных корреспонденциях $\{d_{ij}\}$ (сколько товара перевозится в единицу времени из источника $i$ в сток $j$) идет равновесное формирование способов транспортировки товаров [2]. В результате формируются функции затрат $T_{ij}(\{d_{ij}\})$. Исходя из этих затрат на верхнем уровне, соответствующему медленному времени, решается задача поиска конкурентного равновесия [4, 5] между производителями и потребителями с учетом затрат на транспортировку. В данном случае (см., например, [2]) мы будем иметь дело с адиабатическим исключением быстрых переменных (принцип подчинения Г. Хакена) в случае стохастических динамик. Обоснование имеется в [6].

Различные частные случаи такого рода постановок задач встречались в литературе. Так, например, в классической монографии [7] рассматривается большое количество моделей верхнего уровня, связанных с расчетом матрицы корреспонденций. В монографиях [8 – 10], напротив, внимание сосредоточено на моделях нижнего уровня, в которых с помощью принципа Нэша–Вардропа рассчитывается $T_{ij}(\{d_{ij}\})$. В статье [2] эти модели объединяются для создания единой многоуровневой (в транспортной науке чаще используется термин "многостадийной") модели, включающей в себя и расчет матрицы корреспонденций, и равновесное распределение потоков по путям. В препринте [1] введена терминология, которой мы будем придерживаться и в данной статье, связанная с пунктами производства и потребления, и в отличие от [2, 7] внешняя задача в [1] больше привязана непосредственно к экономике. Но во всех этих случаях можно было обойтись (с некоторыми оговорками в случае [2]) и без понятия конкурентного равновесия, поскольку получающиеся в итоге (популяционные) игры были потенциальными[1], причем имелась и эволюционная интерпретация [11]. Поиск равновесия сводился к решению задачи выпуклой оптимизации, а цены определялись из решения двойственной задачи. В препринте [12] для задачи верхнего уровня была предложена оригинальная конструкция, сводящая поиск конкурентного равновесия к поиску седловой точки (причем, получившаяся игра не была потенциальной в обычном смысле). Тем не менее, в [12] не рассматривалась транспортировка, т.е. не было задачи нижнего уровня.

Целью настоящей работы является предложить такое описание задачи верхнего уровня, включающее в себя описанные выше примеры, которое сводит в итоге поиск конкурентного транспортно-экономического равновесия к поиску седловой точки в выпукло-вогнутой игре. Отметим здесь, что в общем случае поиск конкурентного равновесия сводится к решению задачи дополнительности или (при другой записи) вариационному неравенству [4, 5]. При этом известно, что в общем случае это вычислительно трудные задачи. Однако в ряде случаев экономическая специфика задачи позволяет гарантировать, что полученное вариационное неравенство монотонное. Тогда задача становится существенно привлекательнее в вычислительном плане. В данной статье мы рассматриваем класс задач, в которых вариационные неравенства, возникающие при поиске конкурентных равнове-сий, переписываются в виде седловых задач. Монотонность автоматически обеспечивается правильной выпукло-вогнутой структурой седловой задачи.

---

[1] Вектор-функция затрат, характеризующая затраты при выборе различных стратегий как функция от распределения игроков по стратегиям, является градиентом некоторой скалярной функции.

Опишем вкратце структуру статьи. В п. 2 описывается решение "транспортной" задачи нижнего уровня (ищется равновесное распределение потоков по путям). В п. 3 описывается конструкция равновесного формирования корреспонденций при заданных функциях транспортных затрат. Отметим, что в этих двух пунктах мы фактически работаем только с одним экономическим агентом "Перевозчик" (если смотреть с точки зрения популяционной теории игр, то с агентами одного типа "Перевозчиками"). В п. 3 этот агент(-ы) могут производить товар, неся затраты, и его потреблять, получая выгоду. В п. 4 модель из п. 3 переносится на случай, когда помимо экономического агента "Перевозчик(-и)" в источниках и стоках транспортного графа располагаются независимые от "Перевозчика" новые экономические агенты "Производители" и "Потребители", решающие свои задачи. В заключительном п. 5 модели верхнего и нижнего уровня объединяются в одну общую модель, конкурентное равновесие в которой сводится к поиску седловой точки в выпукло-вогнутой игре.

## 2. Равновесное распределение потоков по путям

Обозначим множество пар $w = (i, j)$ источник-сток $OD$, $x_p$ – поток по пути $p$; $P_w$ – множество путей, отвечающих корреспонденции $w$, $P = \bigcup_{w \in OD} P_w$ – множество всех путей;

$f_e(x) = \sum_{p \in P} \delta_{ep} x_p$ – поток на ребре $e$ (здесь и далее $x = \{x_p\}$, $f = \Theta x$), где $\delta_{ep} = \begin{cases} 1, & e \in p; \\ 0, & e \notin p \end{cases}$;

$\tau_e(f_e)$ – затраты на ребре $e$ ($\tau_e'(f_e) \geq 0$); $G_p(x) = \sum_{e \in E} \tau_e(f_e(x)) \delta_{ep}$ – затраты на пути $p$;

$X = \left\{ x \geq 0: \sum_{p \in P_w} x_p = d_w, w \in OD \right\}$ – множество допустимых распределений потоков по путям, где $d_w$ – корреспонденция, отвечающая паре $w$.

**Определение 1.** *Распределение потоков по путям $x = \{x_p\} \in X$ называется равновесием (Нэша–Вардропа) в популяционной игре $\langle \{x_p\} \in X, \{G_p(x)\} \rangle$, если из $x_p > 0$ ($p \in P_w$) следует $G_p(x) = \min_{q \in P_w} G_q(x)$. Или, что то же самое:*

*для любых $w \in OD$, $p \in P_w$ выполняется $x_p \cdot \left( G_p(x) - \min_{q \in P_w} G_q(x) \right) = 0$.*

**Теорема 1 (см. [2, 8–11]).** *Популяционная игра $\langle \{x_p\} \in X, \{G_p(x)\} \rangle$ является потенциальной. Равновесие $x^*$ в этой игре всегда существует, и находится из решения задачи выпуклой оптимизации*

$$x^* \in \text{Arg}\min_{x \in X} \Psi(f(x)), \qquad (1)$$

*где*

$$\Psi(f(x)) = \sum_{e \in E} \int_0^{f_e(x)} \tau_e(z) dz = \sum_{e \in E} \sigma_e(f_e(x)).$$

Мы оставляем в стороне вопрос единственности равновесия (детали см., например, в [2, 10]). Отметим лишь, что при естественных условиях равновесное распределение потоков по рёбрам $f^*$ единственно. В частности, для этого достаточно, чтобы $\tau_e'(f_e) > 0$ для всех $e \in E$. Если дополнительно $f^* = \Theta x$ однозначно разрешимо относительно $x$ (в реальных транспортных сетях часто случается, что число допустимых для перевозки путей меньше числа рёбер, это как раз и приводит к однозначной разрешимости), то отсюда будет следовать, что равновесное распределение потоков по путям $x^*$ единственно.

Удобно считать [1, 2], что возрастающие функции затрат $\tau_e(f_e) := \tau_e^\mu(f_e)$ зависят от параметра $\mu > 0$, причём

$$\tau_e^\mu(f_e) \xrightarrow[\mu \to 0+]{} \begin{cases} \overline{t_e}, & 0 \leq f_e < \overline{f_e} \\ [\overline{t_e}, \infty), & f_e = \overline{f_e} \end{cases},$$

$$d\tau_e^\mu(f_e)/df_e \xrightarrow[\mu \to 0+]{} 0, \quad 0 \leq f_e < \overline{f_e}.$$

В таком пределе задачу выпуклой оптимизации можно переписать как задачу ЛП [13]:

$$\min_{\substack{f = \Theta x, x \in X \\ f \leq \overline{f}}} \sum_{e \in E} f_e \overline{t_e}.$$

Такого типа транспортные задачи достаточно хорошо изучены [14, 15].

Для дальнейшего будет важно переписать задачу $\min_{x \in X} \Psi(f(x))$ через двойственную [2]:

$$\min_{x \in X} \Psi(f(x)) = \min_{f,x} \left\{ \sum_{e \in E} \sigma_e(f_e) : f = \Theta x, \ x \in X \right\} =$$

$$= \min_{f,x} \left\{ \sum_{e \in E} \max_{t_e \in \text{dom } \sigma_e^*} \left[ f_e t_e - \sigma_e^*(t_e) \right] : f = \Theta x, \ x \in X \right\} =$$

$$= \max_{t \in \text{dom } \sigma^*} \left\{ \min_{f,x} \left[ \sum_{e \in E} f_e t_e : f = \Theta x, \ x \in X \right] - \sum_{e \in E} \sigma_e^*(t_e) \right\} =$$

$$= \max_{t \geq \bar{t}} \left\{ \sum_{w \in OD} d_w T_w(t) - \langle \bar{f}, t - \bar{t} \rangle - \mu \sum_{e \in E} h(t_e - \bar{t}_e, \bar{t}_e, \bar{f}_e, \mu) \right\} =$$

$$\overset{\mu \to 0+}{=} \max_{t \geq \bar{t}} \left\{ \sum_{w \in OD} d_w T_w(t) - \langle \bar{f}, t - \bar{t} \rangle \right\}, \quad (2)$$

где $\sigma_e^*(t_e)$ – сопряженная функция к $\sigma_e(f_e)$, $T_w(t) = \min_{p \in P_w} \sum_{e \in E} \delta_{ep} t_e$ – длина кратчайшего пути из $i$ в $j$ ($w = (i, j)$) на графе, взвешенном согласно вектору $t$, $h(t_e - \bar{t}_e, \bar{t}_e, \bar{f}_e, \mu)$ – сильно выпуклая функция по первому аргументу. При этом

$$\tau_e^\mu(f_e(x(\mu))) \xrightarrow[\mu \to 0+]{} t_e,$$

где $x(\mu)$ – равновесное распределение потоков, рассчитывающееся по формуле (1), а $t = \{t_e\}$ – решение задачи (2), при естественных условиях единственное [2]. Описанный предельный переход позволяет переходить к задачам, в которых вместо функции затрат на ребрах $\tau_e(f_e)$ заданны ограничения на пропускные способности $0 \leq f_e \leq \bar{f}_e$ и затраты $\bar{t}_e$ на прохождения ребер, когда на ребрах нет "пробок" ($f_e < \bar{f}_e$) [2, 13].

Основным для дальнейшего выводом из этого всего является способ (основанный на применении теоремы Демьянова–Данскина [16, 17], как правило, в градиентном варианте в виду единственности $t$) потенциального описания набора $T(d) := \{T_w(t(d))\}$:

$$T(d) = \nabla_d \min_{x \in X(d)} \Psi(f(x)) = \nabla_d \max_{t \geq \bar{t}} \left\{ \sum_{w \in OD} d_w T_w(t) - \langle \bar{f}, t - \bar{t} \rangle - \mu \sum_{e \in E} h(t_e - \bar{t}_e, \bar{t}_e, \bar{f}_e, \mu) \right\}. \quad (3)$$

В [2, 10, 11] приведены эволюционные динамики, приводящие к описанным здесь равновесиям. Отметим, однако, что если рассматривать Logit динамику [2, 11] (ограниченно рациональных агентов с параметром $\tilde{\gamma} > 0$ [18]), то задачу (1) необходимо будет переписать в виде (говорят, что вместо равновесия Нэша–Вардропа ищется стохастическое равновесие [2, 19]):

$$\min_{x \in X} \left\{ \Psi(f(x)) + \tilde{\gamma} \sum_{w \in OD} \sum_{p \in P_w} x_p \ln(x_p / d_w) \right\}. \quad (4)$$

Это замечание понадобится нам в дальнейшем.

В заключение этого раздела отметим, что теорема 1 может быть распространена и на случай, когда затраты на ребрах $\tau_e(f_e; \{f_{\tilde{e}}\})$ удовлетворяют условию потенциальности (частный случай – это когда $\tau_e(f_e; \{f_{\tilde{e}}\}) \equiv \tau_e(f_e)$) [11]:

$$\frac{\partial \tau_e\left(f_e;\{f_{\tilde{e}}\}\right)}{\partial e'} = \frac{\partial \tau_{e'}\left(f_{e'};\{f_{\tilde{e}}\}\right)}{\partial e}.$$

Тогда

$$\Psi(f(x)) = \sum_{e \in E} \int_0^{f_e(x)} \tau_e\left(z;\{\{f_{\tilde{e}}\}_{\tilde{e} \neq e} \cup \{f_e = z\}\}\right) dz.$$

Такого рода обобщение нужно, например, когда пропускные способности узлов зависят от потоков, пересекающих узлы. В случае транспортных потоков такими узлами являются, в частности, перекрестки. Тогда, путем раздутия исходного графа, мы с одной стороны "развязываем узел", сводя затраты прохождения узла по разным путям к затратам прохождения фиктивных (введенных нами) ребер, а с другой стороны приобретаем более общую зависимость $\tau_e\left(f_e;\{f_{\tilde{e}}\}\right)$.

### 3. Равновесный расчет матрицы корреспонденций

В п. 2 матрица корреспонденций $\{d_w\}_{w \in OD}$ была задана по постановке задачи. В данном пункте мы откажемся от этого условия, вводя в источники $O$ производство, а в стоки $D$ – потребление. Агенты "появляются" в тех пунктах производства, произведя товар в которых его можно с выгодой для себя реализовать в каком-нибудь из пунктов потребления. Это означает, что затраты на производство и затраты на транспортировку полностью окупаются последующей выручкой от реализации продукции в пункте потребления. Агенты, которых мы здесь считаем маленькими, будут "приходить" в систему до тех пор, пока существует цепочка (пункт производства–маршрут–пункт потребления), обеспечивающая им положительную прибыль. Важно отметить, что по ходу "наплыва" агентов транспортная сеть становится все более и более загруженной, что может сказываться на затратах на перевозку. В результате прибыль ранее пришедших агентов падает, что побуждает их перераспределяться, т.е. искать более выгодные цепочки. Возникает ряд вопросов. Например, сходится ли такая динамика (точнее семейство динамик, отражающих рациональность агентов) к равновесию? Если сходится, то единственно ли равновесие? Если равновесие единственно, то как его можно эффективно найти (описать)? Попытка ответить на эти вопросы (но, прежде всего, на последний вопрос) для достаточно широкого и важного в приложениях класса задач предпринята в последующей части статьи.

Рассмотрим сначала для большей наглядности отдельно потенциальный случай. А именно тот случай, когда в источнике $i \in O$ производственная функция имеет вид $\sigma_i(f_i)$, где $f_i = \sum_{k:(i,k)=e \in E} f_e = \sum_{j:(i,j) \in OD} d_{ij}$, аналогично для стоков $j \in D$ определим функции полезности со знаком минус $\sigma_j(f_j)$, $f_j = \sum_{k:(k,j)=e \in E} f_e = \sum_{i:(i,j) \in OD} d_{ij}$. Все эти функции считаем выпуклыми. Мы обозначаем эти функции одинаковыми буквами, однако, это не должно вызвать в дальнейшем путаницы в виду характерных нижних индексов. Редуцируем рассматриваемую задачу к задаче п. 2. Рассмотрим новый граф с множеством вершин $O \cup D$, соединенных теми же ребрами, что и в изначальном графе, и с одним дополнительным фиктивным источником и одним дополнительным фиктивным стоком. Этот фиктивный источник соединим со всеми источниками $O$, аналогично фиктивный сток соединим со всеми стоками $D$. Если существует путь из источника $i$ в сток $j$ в исходном графе, то в новом графе прочертим соответствующее ребро с функцией затрат $T_{ij}(d)$. Проведем дополнительное (фиктивное) ребро, соединяющее фиктивный источник с фиктивным стоком, затраты на прохождения которого тождественный ноль. Получим в итоге ориентированный граф путей из источника в сток. Легко понять, что мы оказываемся "почти" в условиях предыдущего пункта (причем с более частным графом – с одним источником и стоком) с точностью до обозначений:

$$\{x_p\} \to \{d_{ij}\}, \{\tau_e(f_e)\} \to \{\sigma'_i(f_i), T_{ij}(d), \sigma'_j(f_j)\}.$$

"Почти", потому что, во-первых, затраты $T_{ij}(d)$ зависят от всего набора $\{d_{ij}\}$, а не только от $d_{ij}$, а во-вторых, не ясно что в данном случае играет роль матрицы корреспонденций (в нашем случае это матрица 1×1, т.е. просто число). Начнем с ответа на второй вопрос. Мы считаем, что в источниках имеется потенциальная возможность производить неограниченное количество продукта, просто в какой-то момент, перестает быть выгодным что-то производить и перевозить. Для этого, собственно, и было введено нулевое ребро, поток по которому обозначим $d_0$. То есть, другими словами, мы должны считать, что $\sum_{(i,j)\in W} d_{ij} + d_0 = \bar{d}$. Если $\bar{d}$ – достаточно большое, то равновесная конфигурация не зависит от того, чему именно равно $\bar{d}$, поскольку не требуется определять $d_0$. С первой проблемой можно разобраться, немного обобщив теорему 1. Предположим, что

$$\exists \ \Phi(d) \text{ - выпуклая}: T(d) = \nabla \Phi(d). \tag{5}$$

Тогда имеет место

**Теорема 2.** *Популяционная игра*

$$\left\langle \{d_{ij}, d_0 \geq 0\}, \{G_{ij}(d) = \sigma'_i(f_i) + T_{ij}(d) + \sigma'_j(f_j), G_0(d) \equiv 0\} \right\rangle,$$

*является потенциальной. Равновесие $d^*$ в этой игре всегда существует (если $\sigma(\cdot)$ – сильно выпуклые функции, то равновесие гарантировано единственно), и находится из решения задачи выпуклой оптимизации*

$$d^* \in \arg\min_{d \geq 0} \tilde{\Psi}(d),$$

$$\tilde{\Psi}\left(d=\{d_{ij}\}\right)=\sum_{i\in O}\sigma_i\left(\sum_{j:(i,j)\in OD}d_{ij}\right)+\sum_{j\in D}\sigma_j\left(\sum_{i:(i,j)\in OD}d_{ij}\right)+\Phi(d). \quad (6)$$

Именно такая конструкция и была рассмотрена в препринте [1] (для многопродуктового рынка). Если искать стохастическое равновесие, то функционал в теореме 2 необходимо энтропийно регуляризовать[2]. Такие конструкции рассматривались (в нерегуляризованном случае), например, в работах [2, 7, 10]. Уже в этих работах можно углядеть необходимость искусственного введения потенциалов (двойственных множителей) в сами функции $\sigma$. А именно, в этих работах предполагается, что все эти функции $\sigma$ – линейные с неизвестными наклонами. Тем не менее, считается, что при этом известно, чему должны равняться в равновесии следующие суммы:

$$\sum_{j:(i,j)\in OD}d_{ij}=L_i, \quad \sum_{i:(i,j)\in OD}d_{ij}=W_j \quad (\sum_{i\in O}L_i=\sum_{j\in D}W_j=N). \quad (7)$$

То есть имеются скрытые от нас (модельера) потенциалы [14] (параметры) $\{\lambda_i^L, \lambda_j^W\}$, которые могут быть рассчитаны исходя из дополнительной информации. Применительно к модели расчета матрицы корреспонденций [2, 7, 10] выписанные дополнительные условия (7) однозначно определяют все неизвестные потенциалы. Однако при этом вместо задачи выпуклой оптимизации мы получаем минимаксную (седловую) задачу выпуклую по $\{d_{ij}\}\geq 0$ и вогнутую, точнее даже линейную, по потенциалам $\{\lambda_i^L, \lambda_j^W\}$:

$$\min_{\substack{\{d_{ij}\}\geq 0 \\ \sum_{(i,j)\in W}d_{ij}=N}}\max_{\{\lambda_i^L,\lambda_j^W\}}\left[\sum_{i\in O}\lambda_i^L\cdot\left(\sum_{j:(i,j)\in OD}d_{ij}-L_i\right)+\sum_{j\in D}\lambda_j^W\cdot\left(W_j-\sum_{i:(i,j)\in OD}d_{ij}\right)+\Phi(d)+\gamma\sum_{(i,j)\in OD}d_{ij}\ln(d_{ij}/N)\right]. \quad (8)$$

Эта задача всегда имеет решение.

### 4. Сетевая модель алгоритмического рыночного поведения

В данном пункте мы предложим сетевой вариант модели поиска конкурентного равновесия из препринта [12]. Однако в контексте изложенного в конце прошлого пункта, нам будет удобнее стартовать с двухстадийной модели транспортных потоков [2], приводящей к равновесию, рассчитываемому по формуле (8).

Предположим теперь, что имеется $m$ видов товаров и, дополнительно, имеется $q$ типов материалов (количества которых можно использовать в единицу времени ограниче-ны вектором $b$), использующихся в производстве. В источниках располагаются производители товаров, а в стоках потребители. Мы считаем, что любой производитель товара, одновременно, является и потребителем, т.е. $O \subseteq D$. Обозначим через $y$ – вектор цен (руб.) на материалы; $\lambda_i^L$ – вектор цен (руб.), по которым производитель продает товары перевозчику в пункте производства $i$, $\lambda_j^W$ – вектор цен (руб.), по которым потребитель покупает товары у перевозчика в пункте потребления $j$. Опишем каждого экономического агента:

---
[2] К сожалению, строгое обоснование (теорема 11.5.12 [11]) имеется только в случае известного (фиксированного) значения $\bar{d}$ (при этом можно считать $d_0=0$).

### *i*-й Производитель

$U_i \subset \mathbb{R}_+^m$ – максимальное технологическое множество (замкнутое, выпуклое);

$\alpha_i \in [0,1]$ – уровень участия;

$\chi_i(\alpha_i U_i) = \alpha_i \chi_i(U_i)$ – постоянные технологические производственные затраты (руб.) при уровне участия $\alpha_i$ (в единицу времени);

$L_i \in \alpha_i U_i$, $[L_i]_k$ – количество произведенного $k$-го продукта (в единицу времени);

$A_i$, $[A_i]_{kl}$ – количество затраченного $l$-го продукта при производстве единицы $k$-го продукта;

$c_i$, $[c_i]_k$ – затраты (руб.) на производство единицы $k$-го продукта;

$R_i$, $[R_i]_{kl}$ – количество затраченного $k$-го материала для приготовления единицы $l$-го продукта.

Описанный "Производитель" решает задачу:

$$\max_{\substack{L_i \in \alpha_i U_i \\ \alpha_i \in [0,1]}} \left\{ \langle \lambda_i^L, L_i \rangle - \chi_i(\alpha_i U_i) - \langle \lambda_i^W, A_i L_i \rangle - \langle c_i, L_i \rangle - \langle y, R_i L_i \rangle \right\} =$$

$$= \max_{L_i \in U_i} \left\{ \left( \langle \lambda_i^L - c^i - A_i^T \lambda_i^W - R_i^T y, L_i \rangle - \chi_i(U_i) \right)_+ \right\}.$$

### *j*-й Потребитель

Предположим, что каждый продукт имеет $s$ различных свойств (своеобразных полезностей). Это может быть, например, содержание витаминов, белков, жиров, углеводов.

$Q_j$, $[Q_j]_{kl}$ – вклад единицы $l$-го продукта в удовлетворение $k$-го свойства;

$\sigma_j$, $[\sigma_j]_k$ – минимально допустимый уровень удовлетворения $k$-го свойства (в единицу времени);

$\beta_j \in [0,1]$ – уровень участия;

$V_j = \{W_j \in \mathbb{R}_+^m : Q_j W_j \geq \sigma_j\}$ – допустимое множество при полном участии;

$W_j \in \beta_j V_j$, $[W_j]_k$ – количество потребленного $k$-го продукта (в единицу времени);

$\tau_j$ – постоянный доход (руб.) при полном участии (в единицу времени).

Описанный "Потребитель" решает задачу:

$$\max_{\substack{W_j \in \beta_j V_j \\ \beta_j \in [0,1]}} \left\{ \beta_j \tau_j - \langle \lambda_j^W, W_j \rangle \right\} = \max_{W_j \in V_j} \left\{ \left( \tau_j - \langle \lambda_j^W, W_j \rangle \right)_+ \right\}.$$

**Перевозчик**

Этот агент решает задачу типа (6), (8):

$$\min_{\{d_{ij}\} \geq 0} \left[ \sum_{i \in O} \left\langle \lambda_i^L, \sum_{j:(i,j) \in OD} d_{ij} \right\rangle - \sum_{j \in D} \left\langle \lambda_j^W, \sum_{i:(i,j) \in OD} d_{ij} \right\rangle + \Phi(d) + \gamma \sum_{(i,j) \in OD} \left( \sum_{k=1}^m [d_{ij}]_k \right) \ln \left( \sum_{k=1}^m [d_{ij}]_k \right) \right],$$

в которой корреспонденции формируются "Производителями" и "Потребителями". Мы считаем, что все товары одинаковы с точки зрения "Перевозчика", т.е. $\Phi(d) := \Phi\left( \left\{ \sum_{k=1}^m [d_{ij}]_k \right\} \right)$ (можно рассматривать и другие варианты). Здесь и далее нам будет удобно писать энтропийную регуляризацию в виде $\gamma \sum_{(i,j) \in OD} \left( \sum_{k=1}^m [d_{ij}]_k \right) \ln \left( \sum_{k=1}^m [d_{ij}]_k \right)$, т.е. опускать $N = \sum_{(i,j) \in OD} \sum_{k=1}^m [d_{ij}]_k$. Точнее полагать $N = 1$ с той же потоковой (физической) размерностью, что и $d$, чтобы под логарифмом была безразмерная величина.[3] При естественных балансовых условиях $\sum_{i \in O} L_i = \sum_{i \in O} A_i L_i + \sum_{j \in D} W_j$ это никак не повлияет на решение задачи.

Проблема здесь в том, что все эти три типа задач завязаны друг на друга посредством векторов цен. Выпишем, как это принято при поиске конкурентных равновесий [4, 5], все имеющиеся **законы Вальраса** (балансовые ограничения + условия дополняющей не-

---

[3] В случае микроскопического обоснования такого рода вариационных принципов (см. п. 5, а также [11]) исходя из рассмотрения соответствующей марковской динамики нащупывания равновесия, мы должны полагать $N \gg 1$, чтобы сделать соответствующий (канонический) скейлинг и перейти к детерминированной постановке.

жесткости), которые накладывают совершенно естественные ограничения на эти векторы цен:

$$\sum_{j:(i,j)\in OD} d_{ij} \leq L_i, \left\langle \lambda_i^L, \sum_{j:(i,j)\in OD} d_{ij} - L_i \right\rangle = 0, \lambda_i^L \geq 0;$$

$$\sum_{k:(k,i)\in OD} d_{ki} \geq W_i + A_i L_i, \left\langle \lambda_i^W, W_i + A_i L_i - \sum_{k:(k,i)\in OD} d_{ki} \right\rangle = 0, \lambda_i^W \geq 0, i \in O;$$

$$\sum_{i:(i,j)\in OD} d_{ij} \geq W_j, \left\langle \lambda_j^W, W_j - \sum_{i:(i,j)\in OD} d_{ij} \right\rangle = 0, \lambda_j^W \geq 0, j \in D \setminus O;$$

$$\sum_{i\in O} R_i L_i \leq b, \left\langle y, b - \sum_{i\in O} R_i L_i \right\rangle = 0, y \geq 0.$$

**Определение 2.** *Набор* $\left\langle \{d_{ij}\}, \{L_i\}, \{W_j\}; y, \{\lambda_i^L\}, \{\lambda_j^W\} \right\rangle$ *называется конкурентным равновесием (Вальраса–Нестерова–Шихмана) если он доставляет решения задачам всех агентов и удовлетворяет законам Вальраса.*

Для того чтобы установить корректность этого определения, подобно [12], введем **условие продуктивности**:

*существуют такие* $\bar{L}_i \in U_i$, $\bar{W}_j \in V_j$, *что* $\sum_{i\in O} \bar{L}_i > \sum_{i\in O} A_i \bar{L}_i + \sum_{j\in D} \bar{W}_j$ *и* $\sum_{i\in O} R_i \bar{L}_i < b$.

**Теорема 3.** *В условиях продуктивности конкурентное равновесие существует и находится из решения правильной выпукло-вогнутой седловой задачи:*

$$\min_{\substack{\{d_{ij}\}\geq 0}} \max_{\substack{\{\lambda_i^L, \lambda_j^W\}\geq 0 \\ y \geq 0}} \min_{\substack{\{L_i \in \alpha_i U_i, \alpha_i \in [0,1]\} \\ \{W_j \in \beta_j V_j, \beta_j \in [0,1]\}}} \left[ \sum_{i\in O} \left( \left\langle \lambda_i^L, \sum_{j:(i,j)\in OD} d_{ij} - L_i \right\rangle + \left\langle \lambda_i^W, \sum_{k:(k,i)\in OD} A_i L_i \right\rangle + \chi_i(\alpha_i U_i) \right) + $$

$$+ \sum_{j\in D} \left( \left\langle \lambda_j^W, W_j - \sum_{i:(i,j)\in OD} d_{ij} \right\rangle - \beta_j \tau_j \right) + \left\langle y, b - \sum_{i\in O} R_i L_i \right\rangle +$$

$$+ \Phi(d) + \gamma \sum_{(i,j)\in OD} \left( \sum_{k=1}^{m} [d_{ij}]_k \right) \ln \left( \sum_{k=1}^{m} [d_{ij}]_k \right) \Bigg] =$$

$$= \max_{\substack{\{\lambda_i^L, \lambda_j^W\}\geq 0 \\ y \geq 0}} \left[ \langle y, b \rangle - \sum_{i\in O} \max_{L_i \in U_i} \left\{ \left( \left\langle \lambda_i^L - c^i - A_i^T \lambda_i^W - R_i^T y, L_i \right\rangle - \chi_i(U_i) \right)_+ \right\} - $$

$$- \sum_{j\in D} \max_{W_j \in V_j} \left\{ \left( \tau_j - \left\langle \lambda_j^W, W_j \right\rangle \right)_+ \right\} + \min_{\{d_{ij}\}\geq 0} \left\{ \sum_{i\in O} \left\langle \lambda_i^L, \sum_{j:(i,j)\in OD} d_{ij} \right\rangle - \sum_{j\in D} \left\langle \lambda_j^W, \sum_{i:(i,j)\in OD} d_{ij} \right\rangle + $$

$$+\Phi(d)+\gamma\sum_{(i,j)\in OD}\left(\sum_{k=1}^{m}[d_{ij}]_k\right)\ln\left(\sum_{k=1}^{m}[d_{ij}]_k\right)\Bigg\}\Bigg]. \qquad (9)$$

Порядок взятия минимума и максимумов можно менять согласно "Sion's minimax theorem" [2, 16, 17, 20].

### 5. Общее конкурентное равновесие

Для того чтобы объединить модели пп. 2, 3, 4 в одну модель рассмотрим формулы (3), (5), (8), (9). Легко понять, что формула (3) как раз и задает тот самый потенциал, существование которого (формула (5)) требуется для справедливости теоремы 2 (неявно это предполагается и в теореме 3), фактически, сводящей поиск конкурентного равновесия к задаче (8), а в общем случае (9).

**Определение 3.** *Набор* $\langle\{x_p\},\{d_{ij}\},\{L_i\},\{W_j\};y,\{\lambda_i^L\},\{\lambda_j^W\}\rangle$ *называется полным (общим) конкурентным равновесием (Вальраса–Нэша–Вардропа–Нестерова–Шихмана) если* $\langle\{d_{ij}\},\{L_i\},\{W_j\};y,\{\lambda_i^L\},\{\lambda_j^W\}\rangle$ *– конкурентное равновесие, а* $\{x_p\}$ *является равновесием (Нэша–Вардропа) при заданном конкурентным равновесием наборе* $\{d_{ij}\}$.

**Теорема 4.** *В условиях продуктивности полное конкурентное равновесие существует и находится из решения правильной выпукло-вогнутой седловой задачи:*

$$\max_{\substack{\{\lambda_i^L,\lambda_j^W\}\geq 0 \\ y\geq 0}}\Bigg[\langle y,b\rangle-\sum_{i\in O}\max_{L_i\in U_i}\left\{\left(\langle\lambda_i^L-c^i-A_i^T\lambda_i^W-R_i^T y,L_i\rangle-\chi_i(U_i)\right)_+\right\}-$$

$$-\sum_{j\in D}\max_{W_j\in V_j}\left\{\left(\tau_j-\langle\lambda_j^W,W_j\rangle\right)_+\right\}+\min_{\{d_{ij}\}\geq 0}\Bigg\{\sum_{i\in O}\left\langle\lambda_i^L,\sum_{j:(i,j)\in OD}d_{ij}\right\rangle-\sum_{j\in D}\left\langle\lambda_j^W,\sum_{i:(i,j)\in OD}d_{ij}\right\rangle+$$

$$+\max_{t\geq\bar{t}}\left\{\sum_{(i,j)\in OD}\left(\sum_{k=1}^{m}[d_{ij}]_k\right)T_{ij}(t)-\langle\bar{f},t-\bar{t}\rangle-\mu\sum_{e\in E}h(t_e-\bar{t}_e,\bar{f}_e,\mu)\right\}+$$

$$+\gamma\sum_{(i,j)\in OD}\left(\sum_{k=1}^{m}[d_{ij}]_k\right)\ln\left(\sum_{k=1}^{m}[d_{ij}]_k\right)\Bigg\}\Bigg]. \qquad (10)$$

Таким образом, поиск общего конкурентного равновесия также сводится к седловой задаче (если мы вынесем все маскимумы и минимумы за квадратную скобку, то получим минимаксную = седловую задачу), имеющей правильную структуру с точки зрения того, что минимум берется по переменным, по которым выражение в квадратных скобках выпукло, а максимум по переменным, по которым выражение вогнуто. Порядок взятия всех максимумов и минимума можно менять согласно "Sion's minimax theorem". В частности, это дает возможность явно выполнить минимизацию по $\{d_{ij}\} \geq 0$, "заплатив" за этого некоторым усложнением получившегося в итоге функционала, который также сохранит правильные выпукло-вогнутые свойства [2].

Мы не будем здесь приводить, что получается после подстановки формулы (3) в формулы (6) и (8). Все выкладки аналогичны, и даже проще. Тем не менее, ссылаясь далее на задачи (6) и (8) мы будем считать, что такая подстановка была сделана.

Такого рода задачи можно эффективно численно решать (причем содержательно интерпретируемым способом), если транспортный граф задачи нижнего уровня (поиска равновесного распределения потоков) не сверх большой [21–24]. Если же этот граф имеет, скажем, порядка $10^5$ ребер, как транспортный граф Москвы и области [10], то требуется разработка новых эффективных методов, учитывающих разреженность задачи и использующих рандомизацию. Мы не будем здесь на этом останавливаться, поскольку планируется посвятить численным методам решения таких задач больших размеров отдельную публикацию. Впрочем, некоторые возможные подходы отчасти освещены в [1, 2]. К сожалению, численный метод, предложенный в [1], не совсем корректен.

Сделаем несколько замечаний в связи с полученным результатом.

Во-первых, в рассмотренных в статье задачах с помощью штрафных механизмов (типа платных дорог) можно добиваться, чтобы возникающие равновесия соответствовали социальному оптимуму. Для этого можно использовать VCG-механизм [25, 26], см. также полную версию статьи [2].

Во-вторых, используя аппарат [1, 2, 12] несложно вычленить из выписанных задач (6), (8), (10) всевозможные цены, тарифы, длины очередей (пробок) – если делаем предельный переход $\mu \to 0+$ и т.п., понимаемые в смысле Л.В. Канторовича, как двойственные множители.

В-третьих, рассматривая два разно масштабных по времени марковских процесса нащупывания равновесной конфигурации можно прийти к решению задач (6), (8) и, с некоторыми оговорками, (10). Например, если в быстром времени динамика перераспределения потоков по путям задается имитационной Logit динамикой [11], а в медленном времени процесс перераспределения корреспонденций (исходя из быстро подстраивающихся затрат $\{T_{ij}(d)\}$) задается просто Logit динамикой [11], то выражение в квадратных скобках (8)) будет играть роль действия в теореме типа Санова, т.е. описывать экспоненциальную концентрацию инвариантной меры марковского процесса с оговоркой, что речь идет о переменных $d$ и $x$ [11]. Аналогичное можно сказать и про $\tilde{\Psi}$ в теореме 2 после подстановки (3). Кроме того, эти же самые функции будут играть роль функций Ляпунова соответствующих прошкалированным (каноническим скейлингом) марковских динамик, приводящих к СОДУ Тихоновского типа [1, 2, 11]. Это также следует из общих результатов работы [27]. Отметим, что относительно нащупывания цен (потенциалов) в задачах (8) и, особенно, (10) работают механизмы похожие на те, которые описаны в классической работе [14]. Другими словами, при фиксированных ценах (потенциалах) динамика соответ-ствует классическим популяционным динамикам нащупывания равновесий [11]. Но из-за того, что потенциалы не известны и, в свою очередь, должны как-то параллельно подбираться предполагается, что в медленном времени экономические агенты переоценивают эти потенциалы исходя из обратной связи (пример имеется в [14]) на то, что они ожидают видеть и то, что они реально видят. Нам не известны строгие результаты, которые бы обосновывали сходимость таких динамик в совокупности. Этому планируется посвятить отдельную публикацию.

В-четвертых, упомянутая выше эволюционная динамика при правильно дискретизации дает разумный численный способ поиска конкурентного равновесия. В частности, упоминаемая имитационная Logit динамика при правильной дискретизации даст метод зеркального спуска / двойственных усреднений, представляющий собой метод проекции градиента с усреднением [23] и без [24], где проекция понимается в смысле "расстояния" Кульбака–Лейблера. Зеркальный спуск можно получить также из дискретизации Logit динамики, если ориентироваться не только на предыдущую итерацию, а на среднее арифметическое всех предыдущих итераций [2, 28]. В работе [28] поясняется некоторая привилегированность этих двух Logit динамик (см. также [2, 11, 18]). Отметим при этом, что Logit динамики может быть проинтерпретирована также как имитационная Logit динамика для потенциальной игры с энтропийно регуляризованным потенциалом [28].

В-пятых, везде выше мы исходили из того, что есть разные масштабы времени. Из-за этого задачи пп. 2 – 4 удалось завязать с помощью формул (3), (5). Однако к аналогичным выводам можно было прийти, если вместо введения разных масштабов времени ввести иерархию в принятии решений [18]. Скажем, сначала пользователь транспортной сети выбирает тип транспорта (личный/общественный), а потом маршрут [2]. Здесь особенно актуальным становятся такие модели дискретного выбора как Nested Logit [18]. А именно, если использовать энтропийную регуляризацию только в одной из этих двух задач разного уровня (иерархии), описанных в пп. 2, 3, то получается обычная (Multinomial) Logit модель выбора (например, в (10) мы регуляризовали только задачу верхнего уровня), но если энтропийно регуляризовать обе задачи, то получится двухуровневая Nested Logit модель выбора [18]. Это означает, что соответствующая Nested Logit динамика в популяционной иерархической игре приводит к равновесию, которое описывается решением задач типа (8), (10) с дополнительной энтропийной регуляризацией задачи нижнего уровня. Несложно показать, что хорошие выпукло-вогнутые свойства задач (8), (10) при этом сохраняются. Да и в вычислительном плане задача не становится принципиально сложнее, особенно если учесть конструкцию "The shortest path problem", описанную в пятой главе монографии [29], см. также [2].

Резюмируем полученные в статье результаты. На конкретных семействах примеров (но, тем не менее, достаточно богатых в смысле встречаемости в приложениях), была продемонстрирована некоторая "алгебра" над различными конструкциями равновесия. Было продемонстрировано, как их можно сочетать друг с другом, чтобы получать все более и более содержательные задачи. Ключевым местом стал переход, связанный с формулой (3), который можно понимать как произведение (суперпозицию) транспортно-экономических моделей, и конструкция задач (8), (9), которую можно понимать как "сумму" моделей. Представляется, что в этом направлении, может возникнуть довольно интересное движение, связанное с вычленением той "минимальной алгебры операций" над моделями, с помощью которой можно было бы описывать большое семейство равновесных конфигураций, встречающихся в различных приложениях.

## Литература